\documentclass{amsart}

\usepackage{amssymb}
\usepackage{amscd}
\usepackage{amsmath}
\usepackage{latexsym}
\usepackage{verbatim}
\usepackage[latin1]{inputenc}
\usepackage[dvips]{graphics}
\relax
\relax
\usepackage[all]{xy}

\newdir{ >}{{}*!/-10pt/\dir{>}}

\newtheorem{Thm}{Theorem}[section]
\newtheorem*{MThm}{Main Theorem}

\theoremstyle{remark}
\newtheorem{Rem}[Thm]{Remark}

\newcommand{\PrfOf}[1]{\noindent\textit{Proof of #1.}}

\newcommand{\bbR}{\mathbb{R}}

\newcommand{\mor}{\mathop{\rm mor}\nolimits}

\newcommand{\res}{\mathop{\rm res}\nolimits}
\newcommand{\Res}{\mathop{\rm Res}\nolimits}
\newcommand{\ind}{\mathop{\rm ind}\nolimits}

\newcommand{\Lind}{\mathop{L\rm ind}\nolimits}

\newcommand{\colim}{\mathop{\rm colim}}
\newcommand{\hocolim}{\mathop{\rm hocolim}}

\newcommand{\nspace}{\hspace*{-.1em}}
\newcommand{\nnspace}{\hspace*{-.05em}}
\newcommand{\nnnspace}{\hspace*{-.02em}}

\newcommand{\cat}[1]{{\mathcal{#1}}}

\newcommand{\Sets}{\mathcal{S}\nspace{\rm e\nnspace t\nnspace s}}
\newcommand{\sSets}{{\rm s}\Sets}
\newcommand{\op}{^{\nnnspace\rm o\nnspace p}}
\newcommand{\Gpds}{\mathcal{G}\nspace{\rm p\nnspace d\nnspace s}}
\newcommand{\Gpdsf}{\Gpds^{\rm f}}
\newcommand{\Cat}{\mathcal{C}\nspace{\rm a}\nnspace{\rm t}}

\newcommand{\Ho}{\mathop{\rm H\nnspace o}\nolimits}
\newcommand{\HsC}[1]{\Ho_{\cat{S}}(\cat{#1})}

\newcommand\UUU[3]{\cat{U}_{#1}(#2,#3)}
\newcommand\USC[2]{\cat{U}_{\cat #1}(\cat #2)}
\newcommand\UC[1]{\cat{U}({\cat #1})}
\newcommand\USCD[3]{\UUU{\cat #1}{\cat #2}{\cat #3}}

\newcommand{\ie}{{\sl i.e.\ }}

\newcommand{\Or}[1]{{\rm Or}(#1)}
\newcommand{\OrGF}[2]{{\rm Or}(#1,#2)}
\newcommand{\Fin}{\mathcal{F\hspace*{-.05em}}in}
\newcommand{\VCy}{\mathcal{V\hspace*{-.05em}C}}

\newcommand{\alg}{{\rm a}\nnspace{\rm l}\nnspace{\rm g}}
\newcommand{\Kalg}{K^{\alg}}
\newcommand{\minfty}{{\scriptscriptstyle{\left<\!-\infty\!\right>}}}
\newcommand{\Lalg}{L^{\minfty}}
\newcommand{\topo}{{\rm t}\nnspace{\rm o}\nnspace{\rm p}}
\newcommand{\Ktop}{K^{\topo}}
\newcommand{\KOtop}{K\!O^{\topo}}

\newcommand{\redRCstar}{{}_{_{\bbR}}\!C^{*}_{\!r}}

\newcommand{\congH}[2]{{}^{#1}\!#2}

\newcommand{\stacksym}[2]{\renewcommand{\arraystretch}{0.1}
\begin{array}{@{}c@{}}
#2\\ #1\\
\end{array}
\renewcommand{\arraystretch}{1}}

\newcommand{\adjtoo}{\ \stacksym{\longleftarrow}{\longrightarrow}\ }

\newcommand{\noloc}{\hspace*{.05em}:\hspace*{-.13em}}

\newcommand{\comma}{\!\searrow\!}

\newcommand{\commaC}[1]{\!\!{\renewcommand{\arraystretch}{0.2}
\begin{array}[b]{c}
{\scriptscriptstyle \;\;\cat #1} \\[-.4em]
\searrow \\
\end{array}
\renewcommand{\arraystretch}{1}}\!\!}

\newcommand{\stackdown}[2]{\renewcommand{\arraystretch}{0.4}
\begin{array}[t]{@{}c@{}}
  #2\\[.06em] {\scriptscriptstyle #1}\\
\end{array}\,
\renewcommand{\arraystretch}{1}}

\newcommand{\dCdot}[2]{\stackdown{\cat{#1}}{#2}}

\begin{document}

%------------------------------------------------------------------------------

\title[Reformulation of the Isomorphism Conjectures]{Model theoretic reformulation
of the Baum-Connes and Farrell-Jones conjectures}

\author{Paul BALMER and Michel MATTHEY}

\address{Department of Mathematics, ETH Zentrum, CH-8092 Zürich, Switzerland}

\email{paul.balmer@math.ethz.ch {\it and} michel.matthey@math.ethz.ch}

\urladdr{http://www.math.ethz.ch/$\sim$balmer {\it and} http://www.math.ethz.ch/$\sim$matthey}

%\address{Department of Mathematics, ETH Zentrum, CH-8092 Zürich, Switzerland}

%\email{}

%\urladdr{www.math.ethz.ch/\~{}matthey}

\thanks{Research supported by Swiss National Science Foundation, grant~620-66065.01}

%------------------------------------------------------------------------------

%\subjclass{Primary ??_??, ??_??; Secondary ??_??, ??_??}

\date{August 6, 2003}

%\keywords{???}

%------------------------------------------------------------------------------

\begin{abstract}
The Isomorphism Conjectures are translated into the language of homotopical algebra, where they resemble Thomason's
descent theorems.
\end{abstract}

%------------------------------------------------------------------------------

\maketitle

%------------------------------------------------------------------------------

\vspace*{-1.5em}
\section{Introduction and statement of the results}

\label{s-Intro-reformulation}%%%%

%------------------------------------------------------------------------------

In \cite{thom}, Thomason establishes that algebraic $K$-theory satisfies Zariski and Nisnevich \emph{descent}. This is
now considered a profound algebraico-geometric property of $K$-theory. In~\cite{bamaI,bamaII}, we have introduced the
sister notion of \emph{codescent}. Here, we prove that each one of the so-called Isomorphism Conjectures (see
\cite{BCH,FaJ3}) among
\begin{itemize}
\item [\textbf{(1)}] the Baum-Connes Conjecture,
\item [\textbf{(2)}] the real Baum-Connes Conjecture,
\item [\textbf{(3)}] the Bost Conjecture,
\item [\textbf{(4)}] the Farrell-Jones Conjecture in $K$-theory,
\item [\textbf{(5)}] the Farrell-Jones Conjecture in $L$-theory,
\end{itemize}
is equivalent to the codescent property for a suitable $K$- or $L$-theory functor.

For a (discrete) group $G$, these conjectures aim at computing, in geometrical and topological
terms, the groups $\Ktop_{*}(C^{*}_{r}G)$, $\KOtop_{*}(\redRCstar G)$, $\Ktop_{*}(\ell^{1}G)$,
$\Kalg_{*}(RG)$ and $\Lalg_{*}(\Lambda G)$ respectively, where $R$ and $\Lambda$ are associative
rings with units, and $\Lambda$ is equipped with an involution. Davis and L\"uck~\cite{dalu}
express these conjectures as follows (the equivalence with the original statements is due to
Hambleton-Pedersen~\cite{hamped}). First, fix one of the Conjectures (1)--(5) and denote by
$K_{*}(G)$ the corresponding $K$- or $L$-group among the five listed above (for (4) and (5),
$R$ and $\Lambda$ are understood). Denote by $\cat{C}:=\Or{G}$ the orbit category of $G$,
whose objects are the quotients $G/H$ with $H$ running among the subgroups of $G$, and the morphisms
are the left-$G$-maps. Let $\cat{D}:=\OrGF{G}{\VCy}$ be the full subcategory of $\Or{G}$ on those
objects $G/H$ for which $H$ is virtually cyclic. We sometimes write $\cat{C}_{G}$ and $\cat{D}_{G}$
to stress the dependence on the group $G$. Then, a suitable functor $X_{G}\colon\cat{C}\longrightarrow
\cat{S}$ is constructed, where $\cat{S}$ denotes the usual stable model category of spectra (of compactly
generated Hausdorff spaces), for which the weak equivalences are the stable ones. This functor $X_{G}$
has the property that $\pi_{*}\big(X_{G}(G/H)\big)$ is canonically isomorphic to $K_{*}(H)$ for all
$H\leqslant G$. Then, the fixed Isomorphism Conjecture for $G$ amounts to the statement that the
following composition, called \emph{assembly map}, is a weak equivalence in $\cat{S}$\,:
$$
\mu^{G}\colon\quad\hocolim_{\cat{D}}\,\res_{\cat{D}}^{\cat{C}}X_{G}\longrightarrow\hocolim_{\cat{C}}
\,X_{G}\stackrel{\sim}{\longrightarrow}\colim_{\cat{C}}\,X_{G}\stackrel{\cong}{\longrightarrow}X_{G}(G/G)\,.
$$
We turn to homotopical algebra. First, we denote by $\USCD{S}{C}{D}$ the model category
$\cat{S}^{\cat{C}}$ of functors $\cat{C}\to\cat{S}$, where the weak equivalences and fibrations
are the \emph{$\cat{D}$-weak equivalences} and \emph{$\cat{D}$-fibrations} respectively, \ie they
are defined $\cat{D}$-objectwise. See details in~\cite[\S\,\ref{s-UCD}]{bamaI},
for instance. For a diagram $X\in\cat{S}^{\cat{C}}$, we let $\xi_{X}\colon QX\longrightarrow X$ be
the cofibrant replacement of $X$ in $\USCD{S}{C}{D}$. As in~\cite[\S\,\ref{s-cod}]{bamaI}, we say
that \emph{$X$ satisfies $\cat{D}$-codescent} if the map $\xi_{X}(c)$ is a weak equivalence in $\cat{S}$
for every $c\in\cat{C}$; if this is only fulfilled at some $c_{0}\in\cat{C}$, we say that \emph{$X$ satisfies
$\cat{D}$-codescent at $c_{0}$}. For a conceptual approach to codescent and a parallel with descent,
see~\cite[\S\S~\ref{s-IntroI} and~\ref{s-CoDvsD}]{bamaI}.
Let $\USC{S}{C}$ be the model structure on $\cat{S}^{\cat{C}}$ with the
$\cat{C}$-weak equivalences and $\cat{C}$-fibrations; we define $\USC{S}{D}$ on $\cat{S}^{\cat{D}}$
similarly. We denote by $\HsC{C}$ and $\HsC{D}$ the homotopy category of $\USC{S}{C}$
and $\USC{S}{D}$ respectively. As in \cite[Prop.\ \ref{Ho-incl-prop}]{bamaI}, we have the derived
adjunction of the Quillen adjunction $\ind_{\cat{D}}^{\cat{C}}\colon\UC{D}\adjtoo\UC{C}\noloc
\res_{\cat{D}}^{\cat{C}}$, namely
$$
\Lind_{\cat{D}}^{\cat{C}}\colon\HsC{D}\adjtoo\HsC{C}\noloc\Res_{\cat{D}}^{\cat{C}}\,.
$$

For the sequel, fix a group $G$ and one of the Isomorphism Conjectures (1)--(5); let
$X_{G}\in\cat{S}^{\cat{C}}$ be the corresponding functor. Keep the other notations as above.

\begin{Thm}
\label{main1-thm}%%%%
The following statements are equivalent\,:
\begin{itemize}
\item [(i)] $G$ satisfies the considered Isomorphism Conjecture;
\item [(ii)] the corresponding functor $X_{G}\in\USCD{S}{C}{D}$ satisfies $\cat{D}$-codescent at
$G/G\in\cat{C}$.
\end{itemize}
\end{Thm}

\begin{Thm}
\label{main3-thm}%%%%
For subgroups $L\leqslant H\leqslant G$, the following statements
are equivalent\,:
\begin{itemize}
\item [(i)] $X_{H}\in\USCD{S}{C_{H}}{D_{H}}$ satisfies $\cat{D}_{H}$-codescent at $H/L\in\cat{C}_{H}$;
\item [(ii)] $X_{G}\in\USCD{S}{C_{G}}{D_{G}}$ satisfies $\cat{D}_{G}$-codescent at $G/L\in\cat{C}_{G}$.
\end{itemize}
\end{Thm}

In fact, by general results of \cite{bamaI} (without invoking \ref{main1-thm} above),
if $X_{G}$ satisfies $\cat{D}_{G}$-codescent, then $X_{H}$ satisfies
$\cat{D}_{H}$-codescent for every subgroup $H\leqslant G$.

\begin{MThm}
\label{main2-thm}%%%%
The following statements are equivalent\,:
\begin{itemize}
\item [(i)] every subgroup $H$ of $G$ satisfies the considered Isomorphism Conjecture;
\item [(ii)] the corresponding functor $X_{G}\in\USCD{S}{C}{D}$ satisfies $\cat{D}$-codescent;
\item [(iii)] up to isomorphism, the image of $X_{G}$ in $\HsC{C}$ belongs to
$\Lind_{\cat{D}}^{\cat{C}}\big(\HsC{D}\big)$.
\end{itemize}
\end{MThm}

Note that the usual Baum-Connes and Bost Conjectures are stated with \emph{finite} subgroups
instead of virtually cyclic ones, but this is known to be equivalent. So, in these cases, we could as well set $\cat{D}_{G}:=\OrGF{G}{\Fin}$ instead of $\OrGF{G}{\VCy}$.

\begin{Rem}
\label{any-cof-repl-rem}%%%%
Let $X\in\cat{S}^{\cat{C}}$ be a diagram and let $\zeta_{X}\colon\mathcal{Q}X\longrightarrow X$
be an arbitrary \emph{cofibrant approximation} of $X$ in $\USCD{S}{C}{D}$, namely, $\zeta_{X}$
is merely a $\cat{D}$-weak equivalence and $\mathcal{Q}X$ is cofibrant in $\USCD{S}{C}{D}$. Then,
$X$ satisfies $\cat{D}$-codescent at some object $c\in\cat{C}$ if and only if $\zeta_{X}(c)$
is a weak equivalence in $\cat{S}$, see \cite[Prop.\ 6.5]{bamaI}. This illustrates the flexibility
of the codescent-type reformulation of the Isomorphism Conjectures, namely, every such cofibrant
approximation of $X_{G}$ yields a possibly very different assembly map that can be used to test
the considered
conjecture.
\end{Rem}

%------------------------------------------------------------------------------

\section{The proofs}

\label{s-prf-main-results}%%%%

%------------------------------------------------------------------------------

%
Let $\Gpdsf$ be the category of groupoids with \emph{faithful} functors. For the considered conjecture,
by \cite{dalu,joa}, there exists a \emph{homotopy functor} $\mathcal{X}\colon\Gpdsf\longrightarrow\cat{S}$,
\ie $\mathcal{X}$ takes equivalences of groupoids to weak equivalences, such that $X_{G}$ is the composite
$$
X_{G}\colon\quad\cat{C}=\Or{G}\stackrel{\iota}{\longrightarrow}\Gpdsf\stackrel{\mathcal{X}}{\longrightarrow}
\cat{S}.
$$
The functor $\iota$ takes $G/H$ to its \emph{$G$-transport groupoid} $\overline{G/H}^{G}$ with the set
$G/H$ as objects and with $\{g\in G\,|\,gg_{1}H=g_{2}H\,\}$ as morphisms from $g_{1}H$ to $g_{2}H$.
Moreover, the functor $\mathcal{X}$ takes values in cofibrant spectra, so that $X_{G}$ is
$\cat{C}$-objectwise cofibrant.

\smallskip

Let $\Cat$ be the category of small categories and $\sSets$ that of simplicial
sets. Denote by $\otimes_{\cat{D}}\colon\sSets^{\cat{D}\op}\!\times\cat{S}^{\cat{D}}
\longrightarrow\cat{S}$ the tensor product over $\cat{D}$ induced by the simplicial
model structure on $\cat{S}$, where $K\in\sSets$ ``acts'' on $E\in\cat{S}$ by $|K|_{+}
\wedge E$.

\smallskip

\PrfOf{Theorem \ref{main1-thm}}
{\sl A priori}, to test whether $X_{G}$ satisfies $\cat{D}$-codescent at some $c\in\cat{C}$ requires
a thorough understanding of the usually mysterious cofibrant \emph{replacement} of $X_{G}$.
A key point here is the freedom to use \emph{any} cofibrant \emph{approximation} instead, see
Remark~\ref{any-cof-repl-rem}. We provide in \cite[\S\,\ref{s-var-cof-approx}]{bamaII}
a general construction of cofibrant approximations in $\USCD{S}{C}{D}$, one of which is exactly
suited for our present purposes \cite[Cor.\ \ref{quote-cor}]{bamaII}. Evaluated at the terminal
object $G/G\in\cat{C}$, this cofibrant approximation $\zeta_{X_{\nnspace G}}\colon\mathcal{Q}X_{G}
\longrightarrow X_{G}$ is a certain map (described at the end of the proof)
$$
\zeta_{X_{\nnspace G}}(G/G)\colon\quad\mathcal{Q}X_{G}(G/G)=B(?\comma\cat{D})\op\dCdot{?\in D}{\otimes}
\res_{\cat{D}}^{\cat{C}}X_{G}(?)\longrightarrow X_{G}(G/G)\,.
$$
Indeed, using the notations of \cite[Not.\ \ref{commaC-not}]{bamaII}, this follows from
the canonical identification $(?\comma\cat{D}\commaC{C}G/G)\op=(?\comma\cat{D})\op$ of diagrams in
$\Cat^{\cat{D}\op}$ and from the fact that $X_{G}$ is $\cat{C}$-objectwise cofibrant. By definition
of the homotopy colimit, we have
$$
B(?\comma\cat{D})\op\dCdot{?\in D}{\otimes}\res_{\cat{D}}^{\cat{C}}X_{G}(?)=\hocolim_{\cat{D}}\,\res_{\cat
{D}}^{\cat{C}}X_{G}\,.
$$
So, it suffices to show that $\zeta_{X_{\nnspace G}}(G/G)$ coincides with the assembly map
$\mu^{G}$. In the notations of \cite[Not.\ \ref{mor-eta-not}]{bamaII}, we have $\underline
{\mor_{\cat{D},\cat{C}}(?,G/G)}=*$ in $\sSets^{\cat{D}\op}$ (the constant diagram with value the point). By \cite[Lem.\ \ref{lambda-iso-lem}]{bamaII}, the spectrum ${*}\,\otimes_{\cat{D}}\res_{\cat{D}}^{\cat{C}}X_{G}$ identifies with
$\ind_{\cat{D}}^{\cat{C}}\res_{\cat{D}}^{\cat{C}}X_{G}(G/G)$. Letting $\epsilon$ denote
the counit of the adjunction $(\ind_{\cat{D}}^{\cat{C}},\res_{\cat{D}}^{\cat{C}})$, it is
routine to verify that there is a canonical commutative diagram
$$
\xymatrix @C=2em @R=1.7em{
\hocolim_{\cat{D}}\res_{\cat{D}}^{\cat{C}}X_{G} \ar@{=}[rr] \ar@{=}[d] & & B(?\comma\cat{D})\op
\otimes_{\cat{D}}\res_{\cat{D}}^{\cat{C}}X_{G}(?) \ar[d] \\
\hocolim_{\cat{D}}\res_{\cat{D}}^{\cat{C}}X_{G} \ar[r] \ar[d] & \colim_{\cat{D}}\res_{\cat{D}}^
{\cat{C}} X_{G} \ar[r]^{\cong} \ar[d] & {*}\,\otimes_{\cat{D}}\res_{\cat{D}}^{\cat{C}}X_{G}
\ar[d]^{\epsilon_{X_{\nnspace G}}(G/G)} \\
\hocolim_{\cat{C}}X_{G} \ar[r]^-{\sim} & \colim_{\cat{C}}X_{G} \ar[r]^-{\cong} & X_{G}(G/G) \\
}
$$
The composition of the first column followed by the last row is the assembly map $\mu^{G}$.
The composition in the last column is $\zeta_{X_{\nnspace G}}(G/G)$, see \cite[Cor.\ \ref{quote-cor}]{bamaII}.
\qed

\smallskip

More generally, one can prove that the ``\,$(X,\mathcal{F},G)$-Isomorphism Conjecture\,''
of \cite[Def.\ 5.1]{dalu} is equivalent to $X$ satisfying $\OrGF{G}{\mathcal{F}}$-codescent
at $G/G$, for any objectwise cofibrant diagram $X\in\cat{S}^{\Or{G}}$ and any family $\mathcal{F}$
of subgroups of $G$.

\smallskip

For $g\in G$ and $H\leqslant G$, we write $\congH{g}{H}:=gHg^{-1}$. In the orbit category
$\Or{G}=\cat{C}_{G}$, for an element $g\in G$ such that $\congH{g}{H}\leqslant K$ for some
subgroups $H$ and $K$ of $G$, we designate by the right coset $Kg$ the morphism $G/H\longrightarrow
G/K$ taking $\tilde{g}H$ to $\tilde{g}g^{-1}K$.

\medskip

\PrfOf{Theorem \ref{main3-thm}}
Consider the functor $\Phi\colon\cat{C}_{H}\longrightarrow\cat{C}_{G}$ taking a coset $H/L\in\cat{C}_{H}$ to
$G/L$. For any $L\leqslant H$, we have canonical equivalences of groupoids in $\Gpdsf$
$$
\overline{H/L}^{H}\stackrel{\sim}{\longleftarrow}\overline{L}\stackrel{\sim}{\longrightarrow}
\overline{G/L}^{G}\,,
$$
where $\overline{L}$ is $L$ viewed as a one-object groupoid. Since $\mathcal{X}$ is a homotopy functor,
one checks that there is a canonical zig-zag of two $\cat{C}_{H}$-weak equivalences between $X_{H}$ and
$\Phi^{*}X_{G}=X_{G}\circ\Phi$ in $\USC{S}{C_{H}}$. By weak invariance of codescent
\cite[Prop.\ \ref{weq-cod-prop}]{bamaI}, $X_{H}$ and $\Phi^{*}X_{G}$ satisfy $\cat{D}_{H}$-codescent
at exactly the same objects $H/L$ of $\cat{C}_{H}$.

Fix an object $H/K\in\cat{D}_{H}$. Let $E_{H/K}\subset G$ be a set of representatives for the quotient
$H\backslash\{g\in G\,|\,\congH{g}{K}\leqslant H\}$. Let $M\gamma\colon\Phi(H/K)=G/K\longrightarrow
G/M=\Phi(H/M)$ be a morphism in $\cat{C}_{G}$ with $M\leqslant H$ (and $\gamma\in G$). It is straightforward
that there is a unique pair $(g,Mh)$ with $g\in E_{H/K}$ and $Mh\in\mor_{\cat{C}_{H}}(H/\congH{g}{K},H/M)$
(namely characterized by $Hg=H\gamma$ and $Mh=M\gamma g^{-1}$) such that $M\gamma$ decomposes in $\cat{C}_{G}$ as
$$
\xymatrix @C=3.5em{
G/K \ar[r]^-{\congH{g}{K}g} \ar@/^18pt/[rr]^-{M\gamma} & G/\congH{g}{K} \ar[r]^-{Mh} & G/M\,.
}
$$
Since $\Phi(\cat{D}_{H})\subset\cat{D}_{G}$, this precisely says that
$\Phi$ is a \emph{left glossy morphism of pairs of small categories} in the sense of
\cite[Defs.\ \ref{CDmorph-def} and \ref{l-glossy-def}]{bamaI}. By left glossy invariance
of codescent \cite[Thm.\ \ref{left-glossy-thm}]{bamaI}, $\Phi^{*}X_{G}$
satisfies $\cat{D}_{H}$-codescent at some $H/L\in\cat{C}_{H}$ if and only if $X_{G}$
satisfies $\cat{D}_{1}$-codescent at $G/L\in\cat{C}_{G}$, where $\cat{D}_{1}:=\Phi(\cat{D}_{H})$.
Set $\cat{D}_{2}:=\cat{D}_{G}$ and fix $H/L\in\cat{C}_{H}$. For $i=1,2$, consider the full
subcategory $\cat{E}_{i}$ of $\cat{D}_{i}$ given by
$$
\cat{E}_{i}:=\big\{G/K\in\cat{D}_{i}\,\big|\,\mor_{\cat{C}_{G}}(G/K,G/L)\neq\varnothing\big\}\,.
$$
By the Pruning Lemma \cite[Thm.\ \ref{prun-obj-thm}]{bamaI},
$X_{G}$ satisfies $\cat{D}_{1}$-codescent at $G/L$ if and only if it satisfies
$\cat{E}_{1}$-codescent at $G/L$. Since $L\leqslant H$, every object of $\cat{E}_{1}$
is isomorphic, inside $\cat{C}_{G}$, to some object of $\cat{E}_{2}$ and conversely;
in other words, $\cat{E}_{1}$ and $\cat{E}_{2}$ are \emph{essentially equivalent}
in $\cat{C}_{G}$, in the sense of \cite[Def.\ \ref{essiso-def}]{bamaI}. So, by
\cite[Prop.\ \ref{shakeD-prop}]{bamaI}, $X_{G}$ satisfies $\cat{E}_{1}$-codescent
at $G/L$ if and only if it satisfies $\cat{E}_{2}$-codescent at $G/L$. By the Pruning
Lemma again, $X_{G}$ satisfies $\cat{E}_{2}$-codescent at $G/L$ if and only if it satisfies
$\cat{D}_{2}$-codescent at $G/L$, \ie $\cat{D}_{G}$-codescent at $G/L$.

In total, we have proven that $X_{H}$ satisfies $\cat{D}_{H}$-codescent at an object
$H/L\in\cat{C}_{H}$ if and only if $X_{G}$ satisfies $\cat{D}_{G}$-codescent at $G/L$,
as was to be shown.
\qed

\smallskip

\PrfOf{the Main Theorem}
The equivalence between (i) and (ii) follows from Theorems \ref{main1-thm} and \ref{main3-thm};
(ii) and (iii) are equivalent by \cite[Thm.\ \ref{Lind-res-thm}]{bamaI}.
\qed

%------------------------------------------------------------------------------

%------------------------------------------------------------------------------

\end{document}